\documentclass[12pt]{amsart}
\usepackage{amsmath}
\usepackage{amsthm}
\usepackage{amscd}

\newcommand{\Q}{{\mathbb Q}}

\newcommand{\Z}{{\mathbb Z}}
\newcommand{\C}{{\mathbb C}}
\newcommand{\F}{{\mathbb F}}
\newcommand{\A}{{\mathbb A}}

\newcommand{\fg}{{\mathfrak g}}

\newtheorem{theorem}{Theorem}
\newtheorem{lemma}{Lemma}
\newtheorem{corollary}{Corollary}
\newtheorem{proposition}[theorem]{Proposition}
\numberwithin{theorem}{section}
\numberwithin{corollary}{section}
\numberwithin{lemma}{section}

\theoremstyle{definition}

\numberwithin{conj}{section}
\newtheorem{example}{Example}
\newtheorem*{acknowledgement}{Acknowledgement}
\numberwithin{example}{section}

\numberwithin{definition}{section}
\newtheorem{question}{Question}
\numberwithin{question}{section}
\numberwithin{equation}{section}

\theoremstyle{remark}
\newtheorem{remark}{Remark}
\numberwithin{remark}{section}

\begin{document}

\title[Recovering representations]{Recovering modular forms 
and representations from tensor \\and symmetric
 powers}
\date{}
\author{C.~S.~Rajan}

\address{Tata Institute of Fundamental 
Research, Homi Bhabha Road, Bombay - 400 005, INDIA.}
\email{rajan@math.tifr.res.in}


\begin{abstract}
We  consider  the problem of determining the relationship between two
representations knowing that some tensor or symmetric power of the
original represetations coincide. Combined with refinements of strong
multiplicity one, 
we show that if the characters of some  tensor or
symmetric powers of two absolutely irreducible $l$-adic representation
with the algebraic envelope of the image being connected, agree 
at the Frobenius elements corresponding to a set of places of
positive upper density, then the representations are twists of each
other by a finite order character.

\end{abstract}

\maketitle

\section{Introduction}
Let $N\geq 1,~k\geq 2$ be positive integers, and $\omega: (\Z/N\Z)^*
\rightarrow \C$, be a character mod $N$, satisfying 
$\omega(-1)=(-1)^k$. Denote by $S(N,k,\omega)^0$
the space of new  forms  on $\Gamma_1(N)$ of weight $k$,
and Nebentypus character $\omega$. 
For $f\in  S_k(N,\omega)^0$ and $p$ coprime to $N$, let
$a_p(f)$ denote the corresponding Hecke eigenvalue of the Hecke
operator at $p$. As an arithmetical 
application  of the results contained in this
paper, we establish the following theorem: 

\begin{theorem} \label{thm:mf}
Let $f\in S_k^n(N,\omega)$ and $f'\in
S_{k'}^n(N',\omega')$ with $k, k'\geq 2$.   
Let $m\geq 1$ be a positive integer, and let 
\[ T:=\{p\mid (p,NN')=1, ~ a_p(f)^m=a_p(f')^m\}.\]
Then the following hold:

a) If $f$ is not a CM-form and  $ud(T)$ is positive,
 then there exists a Dirichlet  character 
$\chi$ of order $m$ such that for all $p$ coprime to $NN'$, we have 
\[ a_p(f)=a_p(f')\chi(p).\]
In particular $k=k'$ and $\omega=\omega'\chi^2$. 

b) Suppose upper  density of $T$ is strictly greater than
 $1-2^{-(2m+1)}$. 
 Then there exists a Dirichlet character $\chi$
such that for all $p$ coprime to $NN'$, we have 
\[ a_p(f)=a_p(f')\chi(p).\]

If we further assume that the conductors $N$ and $N'$ are squarefree
in a) (or in  (b)), then $f=f'$. 

\end{theorem}

When $m=1$, the theorem can be considered as a refinement of strong
multiplicity one and is proved in \cite{rajan98}. 
The theorem was proved for the case $m=2$  by D. Ramakrishnan
\cite{DR}. 

The starting point for the proof of the above theorem are the 
$l$-adic representations $\rho_f: G_{\Q}\to GL_2(\Q_l)$ 
(for suitable rational primes $l$) of the absolute Galois group $G_{\Q}$
of $\Q$  to newforms $f\in S(N, k, \omega)^0$ of
weight $k\geq 2$ by  
the work of Shimura, Igusa, Ihara
and Deligne. This has the property 
that for any rational prime $p$ coprime to $Nl$ the representation
$\rho_f$ is unramified at $p$ and  
\[ {\rm Tr}\rho_f(\sigma_p)= a_p(f),\]
where for a prime $p$ coprime to $N$
$\sigma_p$ denotes the Frobenius conjugacy class at $p$ in the group
$G_{\Q}/{\rm Ker}(\rho_f)$. 
Further, if $\rho_f\simeq \rho_{f'}$ then $f=f'$.  Hence
Theorem \ref{thm:mf} follows provided the
representation $\rho_f$ can be recovered from knowing the
$m^{th}$ tensor product representation 
$\rho_f^{{\otimes}m}:G_{\Q}\to
GL_{2^m}(\Q_l)$.

This leads us to consider the general problem  of recovering $l$-adic
representations from a knowledge of $l$-adic representations
constructed algebraically out of the original representations. 
 Let $F$ be a
non-archimedean local field of characteristic zero. Suppose 
\[ \rho_i:G_K\rightarrow GL_n(F), \quad i=1,2\]
are continuous, semisimple  representations 
of the Galois group $G_K$ into $GL_n(F)$, unramified outside a finite
set $S$ of places containing the archimedean places of $K$. Let $M$ be
an algebraic subgroup of $GL_n$ such that $M(F)$ contains the image
subgroups $\rho_i(G_K), ~i=1,2$. Let 
\[ R: ~M\to GL_m\]
be a rational  representation of $M$ into $GL_m$ defined over $F$.

\begin{question}\label{question} Let $T$  be a  subset  of the set of 
places $\Sigma_K$ of $K$ satisfying,      
\begin{equation}\label{T}
 T=\{v \not\in S \mid{\rm Tr}({R\circ\rho_1}(\sigma_v))= {\rm Tr}({R\circ
\rho_2}(\sigma_v))\}, 
\end{equation}
where $\rho_i(\sigma_v)$ for $v\not\in S$
denotes the Frobenius conjugacy classes lying in
the image. Suppose that $T$ is a  `sufficiently large' set of places of $K$. 
  How are $\rho_1$
and $\rho_2$ related? 

More specifically, under what conditions on $T$,  $R$ or the nature
of the representations $\rho_i$, can we conclude that there exists a
central abelian representation $\chi: G_K\to Z_R(F)$, such that the
representations $\rho_2$ and $\rho_1\otimes \chi$ are conjugate by an
element of $M(F)$?

Further  we would also like to know the answer when we take $R$ to be a
`standard representation', for example if $R$ is taken to be $k^{th}$
tensor, or symmetric or exterior powers of a linear representation of
$M$.  The representation $R\circ\rho$ can be thought of 
as an $l$-adic representation
constructed algebraically from the original representation $\rho$.  

\end{question}

When $M$ is isomorphic to $GL_m$, and $R$ is taken to be the identity
 morphism, then the question is a refinement of strong multiplicity
 one and was considered in  \cite{rajan98}.   In this paper we follow the
 algebraic methods and techniques of our earlier paper
 \cite{rajan98}. 

We  break the general problem outlined above  into two steps.   First
we use the results of \cite{rajan98}, to conclude that $R\circ \rho_1$
and  $R\circ \rho_2$ are isomorphic under suitable density hypothesis
on $T$.  We then consider the algebraic envelopes of the $l$-adic
representations, and  we try to answer the question of recovering
rational  representations of reductive group, 
 where the role of the Galois group in the above problem is  replaced
by a  reductive group.

 The latter aspect can be done in a more abstract context. 
Indeed, 
let $\Gamma$ be an abstract group,  $F$ be an arbitrary field of characteristic
zero,  and let $\rho_1, ~\rho_2:\Gamma\to GL_n(F)$ be
representations of $\Gamma$. Let 
\[ \rho:=\rho_1\times \rho_2:\Gamma\to GL_n(F)\times GL_n(F),\]
be the product representation. Define $G$ (resp. $G_i, ~i=1,2$) as the
algebraic envelope (equivalently the Zariski closure) of the image
group $\rho(\Gamma)$ (resp. $\rho_i(\Gamma)$) in the algebraic group
$GL_n\times GL_n$ (resp. $GL_n$). We continue to   denote by $\rho_i$
the two projection morphisms from $G\to G_i$.  The basic idea of  our
approach is to replace $\Gamma$ by the algebraic envelope  of
$\rho(\Gamma)$ and to  recover the representations of the algebraic
group $G$.  This allows us to apply algebraic, in particular 
 reductive group theoretic
techniques towards  a solution of the problem. For example, the
algebraic methods allow us to work over complex numbers, and also to
work with compact forms of the reductive group. The
theorems  are when $G$ is connected; this has the
effect of replacing $\Gamma$ by a subgroup $\Gamma'$ of finite index in
$\Gamma$.  To go from $\Gamma'$ to $\Gamma$, we need some assumptions
on the nature of the representation restricted to $G_1^0$ (see Lemma
\ref{lem:extension}.

For the arithmetical
applications, typically further information on the nature of the
$l$-adic representation will be required. For the application to
modular forms, we require that the associated $l$-adic representation
$\rho_f$ is semisimple; the algebraic envelope of the image contains a
maximal torus when the weight is at least two; when the form is
non-CM, then the algebraic envelope is the full $GL_2$. 

The contents of the paper are as follows: we first consider the cases
when  $R$ is either the  symmetric, tensor power, adjoint and
twisted tensor product (Asai) 
representations of the ambient group $GL_n$.  We also discuss the
situation when the representation is absolutely irreducible and the
algebraic envelope of the image is not connected, allowing us to treat
CM forms. 

Specialising to  modular forms we generalise the results of
 Ramakrishnan \cite{DR},   where we
consider arbitrary $k^{th}$ powers of the eigenvalues for a natural
number $k$, and also $k^{th}$ symmetric powers in the eigenvalues of
the modular forms, and the Asai representations.

\begin{acknowledgement} The initial idea for this paper was conceived
when the author was visiting Centre de Recherches Math\'{e}matiques,
Montr\'{e}al during the Special Year on Arithmetic held in 1998, and
my sincere  thanks to  CRM for their  hospitality and support. I thank
M. Ram Murty for useful discussions and 
the  invitation to visit CRM, and to Gopal Prasad for
 the reference to the work   of Fong and
Greiss. Some of these results were indicated in \cite{Ra2}.
 The arithmetical 
application of  Theorem \ref{thm:special} d),
to generalized Asai representations was suggested by D. Ramakrishnan's
work, who had earlier proved a similar result for the usual degree two
Asai representations, and I thank him for conveying to me his
results. 
\end{acknowledgement}

\section{Tensor and Symmetric powers, Adjoint and  Asai 
 representations}
We recall the basic setup:  $\Gamma$ is an abstract group,  $F$ is  a
 field of characteristic
zero,   $\rho_1, ~\rho_2:\Gamma\to GL_n(F)$ are
two representations of $\Gamma$,  
$ \rho:=\rho_1\times \rho_2:\Gamma\to GL_n(F)\times GL_n(F)$
is  the product representation, and $R:GL_n\to GL_m$ a rational
representation with kernel contained in the centre of $GL_n$. 
  Define $G$ (resp. $G_i, ~i=1,2$) as the
algebraic envelope (equivalently the Zariski closure) of the image
group $\rho(\Gamma)$ (resp. $\rho_i(\Gamma)$) in the algebraic group
$GL_n\times GL_n$ (resp. $GL_n$). For an algebraic group $G$, $G^0$
will denote the connected component of the identity in $G$. 

 We now specialize $R$ to some familiar
representations. 
For a linear representation $\rho$ of a group $G$ into
$GL_n$, let $T^k(\rho), ~S^k(\rho),
~E^k(\rho)~(k\leq n), ~{\rm Ad}(\rho)$ be respectively the $k^{th}$ tensor,
symmetric, exterior product and adjoint  representations of $G$. 

\subsection{Tensor powers}
\begin{proposition} \label{prop:tensorconn}
Let $G$ be a connected algebraic group over a
characteristic zero base field $F$, and let $\rho_1, ~\rho_2$ be
finite dimensional semisimple representations of $G$ into $GL_n$. Suppose that 
\[ T^k(\rho_1)\simeq T^k(\rho_2)\]
for some $k\geq 1$. Then $\rho_1 \simeq \rho_2$. 
\end{proposition}

\begin{proof} We can work over $\C$. Let $\chi_{\rho_1}$ and
$\chi_{\rho_2}$ denote respectively the characters of $\rho_1$ and
$\rho_2$.  Since $\chi_{\rho_1}^k=\chi_{\rho_2}^k$, the characters 
$\chi_{\rho_1}$
and $\chi_{\rho_2}$ differ by a $k^{th}$ root of unity. Choose a
connected neighbourhood $U$ of the identity in $G(\C)$, where the
characters are non-vanishing. Since $\chi_{\rho_1}(1)= \chi_{\rho_2}(1)$,
and the characters differ by a root of unity, we have $\chi_{\rho_1}=
\chi_{\rho_2}$ on $U$. Since they are rational functions on $G(\C)$
and $U$ is Zariski dense as $G$ is connected, we see that
$\chi_{\rho_1}=\chi_{\rho_2}$ on $G(\C)$. Since the representations
are semisimple, we obtain that $\rho_1$ and $\rho_2$ are equivalent.
\end{proof}

\begin{example}\label{FG}
 The connectedness assumption cannot be dropped. Fong
and Greiss \cite{fong-griess95} (see also Example \ref{blasius}), 
have constructed for infinitely many triples $(n, q, m)$
homomorphisms  of 
$PSL_n(\F_q)$ into $PGL_m(\C)$,  which are elementwise conjugate but
not conjugate as representations. Here $\F_q$ is the finite field with
$q$-elements. Lift  two such homomorphisms  to representations
$$\rho_1, \rho_2 :SL(n, \F_q) \to GL(m,\C).$$
 We obtain  for
each $g\in SL(n,\F_q)$, $\rho_1(g)$ is conjugate to 
$\lambda \rho_2(g)$, with  $\lambda$ a scalar. Let $l$ be an exponent
 of the group $SL(n, \F_q)$. Then $g^l=1$ implies that
 $\lambda^l=1$. Hence 
\[ \chi_{\rho_1}^l=\chi_{\rho_2}^l.\]
Thus the $l^{th}$ tensor powers of $\rho_1$ and $\rho_2$ are
equivalent, but by construction there does not exist a character
$\chi$ of $SL(n,\F_q)$ such that $\rho_2\simeq \rho_1\otimes \chi$.  
\end{example}

\subsection{Symmetric powers}
\begin{proposition} \label{prop:sym}
Let $G$ be a connected reductive algebraic group
over a characteristic zero base field $F$. Let  $\rho_1, \rho_2$ be
finite dimensional representations of $G$ into $GL_n$. Suppose that 
\[ S^k(\rho_1)\simeq S^k(\rho_2)\]
for some $k\geq 1$. Then $\rho_1 \simeq \rho_2$. 
\end{proposition}

\begin{proof} We can work over $\C$. Let $T$ be a maximal torus of
$G$. Since two representations of a reductive group are equivalent if
and only if their collection of weights with respect to a maximal
torus $T$ are the same, it is enough to show that
 the collection of weights of $\rho_1$
and $\rho_2$ with respect to $T$ are the same. By Zariski density or
Weyl's unitary trick, we can work with a compact form of $G(\C)$ with
Lie algebra ${\mathfrak g}_u$. Let ${\mathfrak t}$ be a maximal torus inside
${\mathfrak g}_u$. The weights of the corresponding Lie algebra
representations associated to ${\rho_1}$ and $\rho_2$ are real
valued restricted to ${i\mathfrak t}$. Consequently we can order them with respect to a lexicographic
ordering on the dual of ${i\mathfrak t}$. 

Let $\{\lambda_1,\cdots, \lambda_n\}$ (resp. $\{\mu_1,\cdots, \mu_n\}$) be
the weights of $\rho_1$ (resp. $\rho_2$) with $\lambda_1\geq
\lambda_2\geq \cdots\geq \lambda_n$ (resp. $\mu_1\geq \cdots \geq
\mu_n$). 
The weights of $S^k(\rho_1)$ are composed of elements of the form 
$$\left\{ \sum_{1\leq i\leq n} k_i\lambda_i \!\mid \sum_{1\leq i\leq n}
k_i=k\right\},$$ 
and similarly for $S^k(\rho_2)$.  

By assumption the weights of $S^k(\rho_1)$ and $S^k(\rho_2)$
are same. Since $k\lambda_1$ (resp. $k\mu_1$) is the highest weight of
$S^k(\rho_1)$ (resp. $S^k(\rho_2)$) with respect to the lexicographic
ordering, we have $k\lambda_1=k\mu_1$. Hence $\lambda_1=\mu_1$. 

By induction, assume that for $j<l, ~\lambda_j=\mu_j$. Then the set of
weights $\{\sum_{i<l} k_i\lambda_i\!\mid \sum k_i =k\} $ and 
$\{\sum_{i<l} k_i\mu_i\!\mid \sum k_i =k\} $ are same. Hence the
complementary  sets
$T_1(l)$ (resp. $T_2(l)$) composed of
weights in $S^k(\rho_1)$ (resp. $S^k(\rho_2)$), where at least one
$\lambda_j$ (resp. $\mu_j$) occurs with positive coefficient for some
$j\geq l$ are the same. 

The highest weight in $T_1(l)$ is
$(k-1)\lambda_1+\lambda_l$, and in $T_2(i)$ is $(k-1)\mu_1+\mu_l$. Since
$\lambda_1=\mu_1$, we obtain $\lambda_l=\mu_l$. Hence we have shown
that the collection of weights of $\rho_1$ and $\rho_2$ are the same,
and so the representations are equivalent. 
\end{proof}

\subsection{Adjoint and Generalized Asai representations}
 For a $G$-module $V$,
let ${\rm Ad}(V)$ denote the adjoint $G$-module given by the natural
action of $G$ on ${\rm End}(V)\simeq V^*\otimes V$,  
 where $V^*$ is  
$G$-module dual to  $V$. 

\begin{example} Let $G$ be a semisimple group, and let $V,W$ be non
self-dual irreducible representations of $G$.  Let $V_1=V\otimes W, ~V_2 =V\otimes
W^*$, considered as $G\times G$-modules.  Then as $G\times G$ irreducible modules (or as reducible $G$ modules), ${\rm Ad}(V_1)\simeq {\rm
Ad}(V_2)$, but $V_1$ is neither isomorphic to $V_2$ nor to the dual
$V_2^*$. 
\end{example}

However when $G$ is simple, the following theorem establishing the unique
  factorisation of tensor products of irreducible representations of
  simple Lie algebras is proved in 
  \cite{Ra3}:
\begin{theorem}\label{tensorsimple}
Let $\fg$  be a simple Lie algebra over $\C$. Let $V_1, \cdots, V_n$
and $W_1, \cdots, W_m$ be non-trivial, irreducible, finite dimensional
$\fg$-modules. Assume that there is an isomorphism of the tensor
products, 
\[  V_1\otimes \cdots \otimes V_n\simeq W_1\otimes \cdots \otimes
W_m,\]
as $\fg$-modules. Then $m=n$, and there is a permutation $\tau$ of the
set $\{1, \cdots, n\}$, such that 
\[ V_i\simeq W_{\tau(i)},\]
as $\fg$-modules.

In particular, let $V, ~W$  be irreducible $\fg$-modules and 
\[ {\rm End}(V)\simeq {\rm End}(W),\]
as $\fg$-modules. Then $V$ is either isomorphic to $W$ or  $W^*$. 
\end{theorem}

We  define  now  a generalisation of (pre)-Asai
representations.  Let $\Gamma$ be an abstract group, 
$A$ be a finite set with a map 
\[ \theta: A\to {\rm Aut}(\Gamma).\]
 Given a representation $\rho:\Gamma\to GL_n$, 
define the {\em Asai representation}, 
\[As_{\theta}(\rho)=\otimes_{a\in A}\rho^a,\]
where $\rho^a(\gamma)=\rho(\theta(a)(\gamma)), ~\gamma \in
\Gamma$.      

\begin{example} \label{asai}
Let $K/k$ be a Galois extension with Galois group
$G(K/k)$. Given $\rho$, we can associate the {\em pre-Asai representation}, 
\[As(\rho)=\otimes_{g\in G(K/k)}\rho^g,\]
where $\rho^g(\sigma)=\rho(\tilde{g}\sigma\tilde{g}^{-1})$, where $\sigma \in
G_K$ and $\tilde{g}\in G_k$ is a lift of $g\in G(K/k)$. At an
unramified place $v$ of $K$, which is split completely over a place
$u$ of  $k$, the
Asai character is given by,
\[ \chi_{As(\rho)}(\sigma_v)=\prod_{v|u}\chi_{\rho}(\sigma_v), \]
where $\sigma_v$ is the Frobenius conjugacy class at $v$. 
Hence we get that upto isomorphism, $As(\rho)$ does not
depend on the choice of the lifts $\tilde{g}$. If further $As(\rho)$
is irreducible, and $K/k$ ic cyclic, then $As(\rho)$ extends to a
representation of $G_k$ (called the Asai representation associated to
$\rho$ when $n=2$ and $K/k$ is quadratic).     
\end{example}

\subsection{Exterior powers}
 Let $V$ and $W$ be $G$-modules. It does not seem
possible to conclude in general from an equivalence  of the form
$E^k(V)\simeq E^k(W)$ as $G$-modules, that $V\simeq W$.  Here $E^k(V)$
denotes the exterior $k^{th}$ power of $V$. For example, let $V$ be a
non self-dual $G$-module of even dimension $2n$. If $G$ is semisimple,
then $E^n(V)$ is self-dual, and also that $E^n(V)$ is dual to
$E^n(V^*)$, where $V^*$ denotes the dual of $V$. Hence we have
$E^n(V)\simeq E^n(V^*)$, but $V\not\simeq V^*$. 

{\bf Question:} It would be interesting to know, for the
possible applications in geometry,   the
relationship between two linear representations of  a connected
reductive group $G$ over $\C$ into
$GL_n$, given that their exterior $k^{th}$ power representations for some
positive integer $k<n$ are isomorphic.

\subsection{Extending representations}
The foregoing results allow us to conclude in some cases that the
given representations coincide upon restricting to a subgroup of
finite index in $\Gamma$. Assuming a  suitable `generic' hypothesis on
the nature of the representations, it is possible to determine the
representations on $\Gamma$ itself, but upto twisting by characters: 
\begin{lemma}\label{lem:extension}
 Let  $\rho_1,  ~\rho_2:\Gamma\to GL_n(F)$ be 
representations of $\Gamma$. Suppose there is a normal subgroup  
$\Gamma'$ of $\Gamma$ and a character $\chi':\Gamma'\to F*$
 such that $\rho_2|\Gamma'\simeq \rho_1|\Gamma'\otimes \chi'$.
 Assume further that
$\rho_1|\Gamma'$ is absolutely irreducible, and the character $\chi'$
 is invariant with respect to the action of $\Gamma$ on $\Gamma'$ by
 conjugation.  Then there is a character
$\chi:\Gamma\to F^*$ such that,
\[\rho_2\simeq  \rho_1\otimes \chi.\]
If further all the data are algebraic, then $\chi$ will be an
algebraic character. 
\end{lemma}

\begin{proof} By Schur's lemma, the commuatant of $\rho_1({\Gamma'})$ inside
the algebraic group $GL_n/F$,  is a form of $GL_1$ which over
$\bar{F}$ becomes isomorphic to the group of scalar matrices. 
By Hilbert Theorem 90, the commutant of $\rho_1(\Gamma')$ inside 
$GL_n(F)$ consists of precisely the scalar
matrices. 

For $\sigma\in \Gamma$, let
 $T(\sigma)=\rho_1(\sigma)^{-1}\rho_2(\sigma)\in GL_n(F)$.
 Now for $\tau\in \Gamma'$ and $\sigma\in \Gamma$,
\begin{equation*}
\begin{split}
T(\sigma)^{-1}\rho_1(\tau)T(\sigma) & =
\rho_2(\sigma)^{-1}\rho_1(\sigma)\rho_1(\tau)\rho_1(\sigma)^{-1}\rho_2(\sigma)
\\
&= \rho_2(\sigma)^{-1}\rho_1(\sigma\tau\sigma^{-1})\rho_2(\sigma) \\
& =\rho_2(\sigma)^{-1}\rho_2(\sigma\tau\sigma^{-1})\chi'(\sigma\tau\sigma^{-1})^{-1}\rho_2(\sigma)\\
& = \rho_2(\tau)\chi'(\tau)^{-1}\\
&=\rho_1(\tau)
\end{split}
\end{equation*}
Thus  $T(\sigma)$ is equivariant with respect
to the representation $\rho_1\!\mid_{\Gamma'}$, and hence is given by a scalar
matrix $\chi(\sigma)$. Since $\chi(\sigma)$ is a scalar matrix, it follows
that for $\sigma,~\tau\in \Gamma$,
$\chi(\sigma\tau)=\chi(\sigma)\chi(\tau)$, i.e., $\chi$ is character
of $\Gamma$ into the group of invertible elements $F^*$ of $F$,
and $\rho_2(\sigma)=\chi(\sigma)\rho_1(\sigma)$ for all $\sigma \in
\Gamma$. 

It is clear from the proof that the lemma holds in the algebraic case
too. 
\end{proof}

\begin{remark} The proof indicates that the converse also holds: if
  there exists a character $\chi$, then $\chi'$ is forced to be an
  invariant with respect to the action of $\Gamma$.
\end{remark}

\subsection{} We summarize the results obtained so far when $R$ is
taken to be a special representation.  

\begin{theorem}\label{thm:special}
 Let $\Gamma$ be an abstract group,
  and $F$ be a field of characteristic zero. 
Let $\rho_1,~\rho_2:~\Gamma \to GL_n(F)$ be
semisimple representations of $\Gamma$.  
Let $R$ be a rational  representation of
$GL_n\to GL_m$ such that 
\[R\circ
\rho_1\simeq R\circ\rho_2 .\] 
Then the following holds: 
\begin{enumerate}
\item Suppose $R= T^k~\text{or}~S^k$ for some natural number $k$. 
 Then there  is a subgroup $\Gamma'$ of finite index in $\Gamma$ such that 
\[\rho_1\!\mid \Gamma'\simeq \rho_2 \!\mid \Gamma'.\]
If further $G^0_1$ acts absolutely irreducibly on $F^n$, then 
there is a character
$\chi:\Gamma \to F^*$ such that,
\[\rho_2\simeq \rho_1\otimes \chi.\]

\item Suppose $R=Ad$ be the adjoint representation of $GL_n$. Assume
  further that  the
  Lie algebra of the derived subgroup of $G^0_1$ is simple 
and  acts absolutely irreducibly
  on $F^n$. 
Then there is a character
$\chi:\Gamma \to F^*$ such that,
\[\rho_2\simeq \rho_1\otimes \chi\quad \text{ or} \quad  
\rho_2\simeq \rho_1^*\otimes \chi.\]

\item Let $A$ be a finite set with a map, 
$\theta: A\to {\rm Aut}(\Gamma).$ Let $R=As_{\theta}$. 
Assume
  further that  the
  Lie algebra of the derived subgroup of $G^0_1$ is simple
 and  acts absolutely irreducibly
  on $F^n$.  Then there is a character
$\chi:\Gamma\to F^*$, and an element $a\in A$ such that
\[\rho_2\simeq \rho_1^a\otimes \chi.\]
\end{enumerate}

\end{theorem}

\begin{proof} Let $\rho:=\rho_1\times \rho_2 :\Gamma\to (GL_n\times
GL_n)(F)$ be the product representation. Let
$G=G_{\rho}$. Continue to denote by
$\rho_i$ the two projection morphisms from $G$ to $GL_n$. 

 Applying the above Propositions \ref{prop:tensorconn},
\ref{prop:sym} and \ref{tensorsimple} to $G^0$ and Lemma
\ref{lem:extension} (where
the role of $\Gamma'$ is played by $G_d$ which is a characteristic
subgroup of $G^0$), we obtain  
the theorem for the algebraic group $G$. Upon restricting back to
$\Gamma$ we obtain the theorem for any abstract group.  

\end{proof}

\section{Tensor powers of Induced representations: the nonconnected case}
The foregoing results allow us to recover representations on the full
group $\Gamma$ itself rather than recovering representations upto
restriction to subgroups of finite index, provided that the connected
component $G_1^0$ of the algebric envelope of the image
$\rho_1(\Gamma)$ acts absolutely irreducibly on $F^n$. As Examples
\ref{FG} and \ref{blasius} show,  it is not possible in
general to recover representations upto twisting by characters for the
full group $\Gamma$. In order to get some positive results even when
$G_1^0$ does not act absolutely irreducibly on $F^n$, we restrict to a
specific situation keeping in mind the application to recovering CM
forms from knowing the equality of `tensor powers' of such forms. 
In particular, such representations  
are induced from a subgroup of finite index. 
We first present another application of the algebraic machinery.
\begin{proposition}\label{tensorirr}
 Let $G$ be an algebraic group over a
characteristic zero base field $F$, and let $\rho_1, \rho_2$ be
finite dimensional semisimple representations of $G$ into $GL_n$. 
Suppose that \[ T^k(\rho_1)\simeq T^k(\rho_2)\]
for some $k\geq 1$. Then if  $\rho_1$ is irreducible, so is $\rho_2$. 
\end{proposition}
\begin{proof} The proof follows by base changing to complex numbers
and considering the characters $\chi_1,\chi_2$ 
of the representations restricted to a
maximal compact subgroup. The characters differ by a root of unity in
each connected component, and hence
$<\chi_1,\chi_1>=<\chi_2,\chi_2>$. Hence it follows by Schur
orthogonality relations that if one of the representations is irreducible,
so is the other.
\end{proof}

We continue with the hypothesis of the above proposition.   By Proposition
\ref{prop:tensorconn}, we can assume  
that $\rho_1|G^0$ and  $\rho_2|G^0$ are isomorphic, which we will
denote by $\rho^{0}$. Let 
$\Phi:=G/G^0$ be the group of connected components of $G$. 
 For $\phi\in \Phi$, let $G^{\phi}$ denote the corresponding connected
component (identity component is $G^0$). The groups $G$ and $\Phi$
act on the collection
of representations of $G^0$: 
\[ r^{\phi}(x)=r(g_{\phi}^{-1}xg_{\phi}) \quad x\in G^0, \quad
g_{\phi}\in G^{\phi},\]
where $r$ is a representation of $G^0$. The action does not depend on
the choice of $g_{\phi}\in G^{\phi}$.   Write 
\[\rho^0=\oplus_{i \in I} r_i,\] where $r_i$ are the
representations on the isotypical components, and $I$ is  the
indexing set of the isotypical components of $\rho^{(0)}$. Fix $i_0\in
I$, and  define $G'$
to be the stabilizer of $r_{i_0}$ in $G$. 
With these assumptions we have the following consequence of Clifford theory:
\begin{proposition}\label{prop:clifford} 
 Let $G$ be an algebraic group over a
characteristic zero base field $F$, and let $\rho_1, \rho_2$ be
finite dimensional semisimple representations of $G$ into $GL_n$. 
Suppose that \[ T^k(\rho_1)\simeq T^k(\rho_2)\]
for some $k\geq 1$. Assume further that either $\rho_1$ or $\rho_2$ is
absolutely irreducible. 
 Then
\[\rho_1={\rm Ind}_{G'}^G(r_1')\quad \text{and}\quad  \rho_2={\rm
  Ind}_{G'}^G(r_2') \]
are induced representations respectively from representations
$r_1', ~r_2'$ of  the subgroup $G'$ of $G$. 
\end{proposition} 

Assume now that the constituents $r_i$ of
$\rho^0$ are absolutely irreducible, i.e, the irreducible
representations occur with multiplicity one.
In the notation of the above Proposition \ref{prop:clifford}, the assumption of
multiplicity one on $\rho^0$, implies that $r_2'=r_1'\otimes \chi$ for
some character $\chi\in {\rm Hom}(G', GL_1)$ trivial upon restriction
to $G^0$.
 Assuming the hypothesis of Proposition \ref{tensorirr}, we
would like to know whether $\rho_2$ and $\rho_1$ differ by a
character. This amounts to knowing that the character $\chi$ extends
to a character of $G$, since the representations are induced. 
Assume from now onwards that $G'$ is normal in $G$. Then the question
of extending $\chi$ amounts first to showing that $\chi$ is invariant
and then to show that invariant characters extend.

\begin{remark}
We rephrase the problem in a different language, with the hope that it
may shed further light on the question. For $\sigma\in G$, let
$T(\sigma)=\rho_1(\sigma)^{-1}\rho_2(\sigma)$ be as in Lemma
\ref{lem:extension}. The calculations of Lemma \ref{lem:extension} with now
$\sigma \in G, \tau \in G^0$, show that $T(\sigma)$ takes values in
the commutant $S$ of the image of $\rho^0$, trivial when restricted to
$G^0$.
Since we have assumed that $\rho_1$ is
irreducible, $\Phi=G/G^0$ acts on $S$ via $\sigma_1$ transitively as a
permutation representation  on 
indexing set $I$. Hence $I$ can be taken to be $G/G'$ and 
$S$ is isomorphic to the induced module ${\rm Ind}_{G'}^G(1_{G'})$, where
the $1_{G'}$ denotes the trivial one dimensional representation  of
$G'$.
 We have for $\sigma, \tau \in
G$, 
\[ T(\sigma\tau)=\rho_1(\tau)^{-1}T(\sigma)\rho_1(\tau)T(\tau), \]
i.e, $T$ is a one cocycle on $G$ with values in $S$. 
Since $S$ is induced, the `restriction' map $H^1(G, S)\to H^1(G', F^{*})$
is an isomorphism. The invariants of the $G$-action on $S$ is given by
the diagonal $F^{*}$ sitting inside $S$, and the composite map,
$ {\rm Hom}(G, F^{*}) \to H^1(G, S)\to {\rm Hom}(G', F^{*})$ 
is the restriction
map. To say that $\chi$ extends to a character of $G$, amounts to
knowing that $\chi$ lies in the image of this composite map. 
\end{remark}

\begin{example}\label{blasius}
 We consider the example given by Blasius \cite{blasius94} in
this context. Let $n$ be an odd prime, and let $H_n$ be the  finite
Heisenberg group, with generators $A, ~B, ~C$ subject to  the
relations: $A^n=B^n=C^n=1, ~AC=CA, ~BC=CB, ~AB=CBA$. Let $e_1, \cdots,
e_n$ be a basis for $\C^n$, and let $\xi_n$ be a primitive $n^{th}$
root of unity. For each integer
$a$ coprime to $n$, define the representation $\rho_a:H_n\to GL_n(\C)$
by, 
\begin{equation*}
\begin{split}
 \rho_a(A)e_i& =\xi_n^{(i-1)a}e_i\\
 \rho_a(B)e_i& =e_{i+1}\\
 \rho_a(C)e_i& =\xi_n^{a}e_i,
\end{split}
\end{equation*}
where the notation is that $e_{n+1}=e_1$. It can be seen that $\rho_a$
are irreducible representations, and that the corresponding projective
representations for any pair of integers $a, ~b$ not congruent modulo
$n$,  are inequivalent. Further for any element $h\in H_n$,
the images of $\rho_a(h)$ and $\rho_b(h)$  in $PGL(n, \C)$ are
conjugate. Hence it follows that for some positive integer $k$ (which
we can take to be $n$) the representations  $T^k(\rho_a)$ and
$T^k(\rho_b)$ are isomorphic. 

Let $T$ be the abelian normal
subgroup of index $n$ generated by $A$ and $C$. There exists a
character $\chi_{ab}$ of $T$ such that $\rho_a|T\simeq \rho_b|T\otimes
\chi_{ab}$. From the theory of
induced representations, it can be further checked that $\rho_a|T$ has
multiplicity one, and that $\rho_a$ is induced from a character
$\psi_a$ of $T$. However we have that there does not exist any
character $\eta$ of $H_n$ such that $\rho_a\simeq \rho_b\otimes \eta$.
\end{example}

With this example in mind, we now present a proposition in the positive
direction (which applies in particular to CM forms of weight $\geq 2$): 

\begin{proposition}\label{prop:rec-ind}
With notations and hypothesis as in Proposition  \ref{prop:clifford}, assume
further the following:
\begin{itemize}
\item the subgroup $G'$  is normal in $G$. 
\item the isotypical components of $\rho_1$ restricted to $G^0$ are
  one dimensional. 
\end{itemize} 
 Let $\chi$ be a character of $G'$
such that $r_2'\simeq r_1'\otimes \chi$. 
Then $\chi$ is invariant with respect to the action of $G$ on the
characters of $G'$. In particular, if invariant characters of $G'$
extend to invariant characters of $G$ (which happens if $G/G'$ is
cyclic), then there is a character $\chi$ of $G$ such that,  
\[ \rho_2\simeq \rho_1\otimes \chi.\]
\end{proposition}
\begin{proof}
Restricting $\rho_1$ and $\rho_2$ to $G'$, have by our assumptions
\[ \rho_1|G'=\oplus_{\substack{\phi \in G/G'}}{r_1'}^{\phi},\quad {\rm
and}\quad \rho_2|G'=\oplus_{\substack{\phi \in G/G'}}({r_1'\chi})^{\phi},\]
where $\chi$ is a character of $G'$ trivial on $G^0$, and $r_1',
~r_2'$ are characters of $G'$.  
Let $\tau$ be an element of $G-G'$. 
If $\chi\neq \chi^{\tau}$, choose
an element $\theta$ of $G'$ such that $\chi(\theta)\neq
\chi^{\tau}(\theta)$. Choose an element $g_{\theta}\in G^{\theta}$. 
Since $G^{\theta}$ is connected, we obtain from our assumption 
$T^k(\rho_1)\simeq T^k(\rho_2)$,  that for some $\zeta \in {\mathbf
\mu}_k$, 
\[\zeta\sum_{\phi\in G/G'}(r_1(g_{\theta})r_1(y))^{\phi}=\sum_{\phi\in
G/G'}(r_1(g_{\theta})r_1(y)\chi(g_{\theta}))^{\phi}, \]
for all $y\in G^0$. But the  characters $r_1^{\phi}$ of $G^0$ are linearly
independent, and hence the above equality forces the character $\chi$
to be invariant. 
\end{proof}

\begin{remark} It seems possible that the above proposition holds
  without the second assumption  that $r_1$ be a character. 
\end{remark}

\section{Strong multiplicity one}
Let $K$ be a global field. Denote by  $\Sigma_K$ the set of places of $K$. 
For a nonarchimedean place $v$ of $K$,  let ${\mathfrak p}_v$ denote 
the corresponding prime ideal of ${\mathcal O}_K$,  and
$Nv$ the norm 
of $v$ be the number of elements of the finite field 
${\mathcal O}_K/{\mathfrak p}_v$. The upper density  
$ud(P)$ of a set $P$ of primes of $K$, is defined to be the ratio,
\[
ud(P)={\varlimsup}_{x\rightarrow \infty}\frac
{\#\{v\in \Sigma_K\mid Nv\leq x,~v\in P\}}
{\#\{v\in \Sigma_K\mid Nv\leq x\}}.\]

Suppose  $L$ is  a finite Galois extension of $K$,
with Galois group $G(L/K)$. Let $S$ denote a  subset of $\Sigma_K$, 
containing the archimedean places together with the set of places of $K$
which ramify in $L$. For each place $v$ of $K$ not in $S$, and a place 
$w$ of $L$ lying over $v$, we have a canonical Frobenius element $\sigma_w$
in $G(L/K)$, defined by the following property:
\[ \sigma_w(x)\cong x^{Nv}({\rm mod}\,{\mathfrak p}_w).\]
 The set $\{\sigma_w\mid w|v\}$ form the Frobenius conjugacy 
class in $G(L/K)$, which we continue to denote by $\sigma_v$. 

  Let $M$ denote a connected, reductive 
algebraic group defined over $F$, and let $\rho$ be
 a continous representation of $G_K$ into 
$M(F)$, where $F$ is a non-archimedean 
local field of residue characteristic $l$, where $l\neq {\rm
  char}(K)$.   Let $L$ denote the fixed field
of $\bar{K}$ by the kernel of $\rho.$ Write $L=\cup_{\alpha}L_{\alpha},$ where 
$L_{\alpha}$ are finite extensions of $K$. $\rho$ is said to be unramified 
outside a set of primes $S$ of $K$, if each of the extensions 
$L_{\alpha}$ is an unramified extension of $K$ outside $S$.

We will assume henceforth that all our linear $l$-adic 
representations of $G_K$  are continuous and semisimple, since we need
to determine a linear representation from it's character.  By the
results of \cite{KR}, 
 it follows that the set of ramified primes is of
density zero, and hence arguments involving density as in
\cite{rajan98}, go through essentially unchanged. Hence in what
follows, $S$ will indicate a
a set of primes of density zero, containing the ramified primes of the
(finite) number of l-adic representations under consideration, and the
archimedean places of $K$. 
 
 Let $w$ 
be a valuation on $L$ extending a valuation $v\not\in S.$ The Frobenius 
elements at the various finite layers for the valuation $w\mid_{L_{\alpha}}$ 
patch together to give raise to the Frobenius element $\sigma_w\in G(L/K)$,
and a Frobenius conjugacy class $\sigma_v\in G(L/K).$ Thus $\rho(\sigma_w)$
(resp. $\rho(\sigma_v)$) is a  well defined element (resp. 
conjugacy class) in $M(F)$. 
If $\rho:G_K\to GL_m(F)$ is a linear representation,
let $\chi_{\rho}$ denote the associated character. 
 $\chi_{\rho}({\sigma_v})$ is well
defined for $v$ a place of $K$ not in $S$. 

We recall here the results of \cite{rajan98}. Suppose 
\[ \rho_i:G_K\rightarrow GL_n(F), \quad i=1,2\]
are continuous, semisimple  representations 
of the Galois group $G_K$ into $GL_n(F)$, unramified outside a 
set $S$ (as above) of places containing the archimedean places of $K$. 
 Let $T$  be a  subset  of the set of 
places $\Sigma_K$ of $K$ satisfying,      
\[
 T=\{v \not\in S \mid{\rm Tr}(\rho_1(\sigma_v))= {\rm Tr}(
\rho_2(\sigma_v)), \quad v\not\in S. \}
\]

Consider the following two hypothesis on the upper density of $T$,
depending on the representations $\rho_1$ and $\rho_2$:
\begin{align}
& {\bf DH1:} &  ud(T) &  >1-1/2m^2 \\
& {\bf DH2:} & ud(T) &  >{\rm min}(1-1/c_1,1-1/c_2),
\end{align}
where $c_i=|G_i)/G_i^0|$
 is  the number of connected components of $G_i$.
As a consequence of the refinements of strong multiplicity one proved
in \cite[Theorems 1 and 2]{rajan98}, we obtain
\begin{theorem} \label{thm:smo} 
 i) If $T$ satisfies {\bf DH1},  then $\rho_1\simeq \rho_2$.

ii) If $T$ satisfies {\bf DH2},
then  there is a finite Galois extension $L$ of $K$, such that 
$\rho_1\!\mid_{G_L}\simeq \rho_2\!\mid_{G_L}$.
 The connected component $G_2^0$  is conjugate to $G_1^0$. In
particular if 
either $G_1$ or $G_2$ is connected and $ud(T)$ is positive, then there
is a finite Galois extension $L$ of $K$, such that 
$\rho_1\!\mid_{G_L}\simeq \rho_2\!\mid_{G_L}$.

\end{theorem}
\begin{proof}
For the sake of completeness of exposition, we present a brief outline
of the proof and refer to \cite{rajan98} for more details.
Let $\rho=\rho_1\times \rho_2$. Let $G$ denote the Zariski closure of
the image $\rho_1(G_K)\times \rho_2(G_K)$. Let $G^0$ denote the
connected component of the identity in $G$, and let $\Phi:=G/G^0$
be the group of connected components of $G$. For $\phi\in \Phi$, let
$G^{\phi}$ denote the corresponding connected component. 

Consider the algebraic subscheme 
\[ X=\{(g_1,g_2)\in G\mid {\rm Tr}(g_1)={\rm Tr}(g_2)\}.\]
It is known that if $C$ is a closed, analytic subset of
$G$, stable under conjugation by $G$ and of dimension strictly smaller
than that of $G$, then the set of Frobenius conjugacy classes lying in
$C$ is of density $0$. Using this it follows that the collection of 
Frobenius conjugacy classes lying in
$X$ has a density equal to 
\begin{equation}\label{algdensity}
\lambda= \frac{|\{\phi\in G/G^0\mid G^{\phi}\subset X\}|}{|G/G^0|}.
\end{equation}
Since this last condition is algebraically defined, the above
expression can be calculated after base changing to $\C$. Let $J$
denote a maximal compact subgroup of $G(\C)$, and let $p_1,~p_2$
denote the two natural projections of the product $GL_m\times GL_m$. 
Assume that $p_1$ and $p_2$ give raise to inequivalent representations
of $J$. (i) follows from the inequlities
\begin{equation}\label{orthogonality}
 2\leq \int_J |{\rm Tr}(p_1(j))-{\rm Tr}(p_2(j))|^2d\mu(j)\leq
(1-\lambda)4m^2,
\end{equation}
where $d\mu(j)$ denotes a normalized Haar measure on $J$.
 The first
inequality follows from the orthogonality relations for
 characters. For the second inequality, we observe that the eigenvalues of
 $p_1(j)$ and $p_2(j)$ are roots of
 unity, and hence 
\[ |{\rm Tr}(p_1(j))-{\rm Tr}(p_2(j))|^2\leq 4m^2.\]
Combining this with the expression for the density $\lambda$ given by
equation \ref{algdensity} gives us the second inequality.

To prove (ii), it is enough to show that $H^0\subset X$. Let
$c_1<c_2$. The density hypothesis implies together with the expression
(\ref{algdensity}) for the density, that there is some element of the
form $(1,j)\in J\cap X$. The proof concludes by observing that the
only element in the unitary group $U(m)\subset GL(m,\C)$ with trace
equal to $m$ is the identity matrix, and hence the connected component
of the idenity in $J$ (or $H$) is contained inside $X$. 
\end{proof} 

\begin{remark}
In the automorphic context, assuming the Ramanujan-Petersson
conjectures, it is possible to obtain the inequalities in
(\ref{orthogonality}), by analogous arguments, and thus a proof of
Ramakrishnan's conjecture in the automorphic context. The first
inequality follows from replacing the orthogonality relations for
characters of compact groups, by the Rankin-Selberg convolution of
$L$-functions, and amounts to studying the behvior 
at $s=1$ of the logarithm  of the function, 
\[L(s,`{|\pi_1-\pi_2|^2}'):= \frac{L(s,\pi_1\times
\tilde{\pi}_1)L(s,\pi_2\times \tilde{\pi}_2)}{L(s,\pi_1\times
\tilde{\pi}_2)L(s,\tilde{\pi}_1\times {\pi}_2)},\]
 where $\pi_1$ and $\pi_2$ are unitary,
automorphic representations of $GL_n({\bf A}_K)$ of a number field $K$, and
$\tilde{\pi}_1, ~\tilde{\pi}_2$ are the contragredient representations
of $\pi_1$ and $\pi_2$. The second
inequality follows from the Ramanujan hypothesis. For more details we
refer to \cite{Ra1}. 
 \end{remark}

\section{Applications to Modular Forms}
 As a corollary to Theorem
\ref{thm:special} and Proposition \ref{prop:rec-ind},  we give an
application to  two dimensional representations:

\begin{corollary}  \label{cor:twodim}
Let $K$ be a number field and let $F$ be a non-archimedean local field
of residue characteristic $l$. Let $\rho_1,
~\rho_2:G_K\to GL_2(F)$ be continuous $l$-adic representations as
above. Suppose that the algebraic envelope  $G_1$ of the image  $\rho_1(G_K)$
contains a maximal torus.  Let $R:GL_2\to GL_m$ be either the
symmetric $k^{th}$ power representation $S^k$ or the $k^{th}$ tensor
power representation $T^k$ of $GL(2)$. 
Let $T$  be a  subset  of the set of 
places $\Sigma_K$ of $K$ satisfying,      
\[
 T=\{v \not\in S \mid{\rm Tr}((R\circ \rho_1)(\sigma_v))= {\rm Tr}((R\circ
\rho_2)(\sigma_v)), \quad v\not\in S, \}
\]  
where $S$ is a finite set of places of $K$ containing the archimedean
places of $K$ and the places of $K$ where either $\rho_1$ or $\rho_2$
is ramified.  Then the following holds:
\begin{enumerate}
\item  Assume that the algebraic envelope of the image of $\rho_1$
contains $SL_2$. If $ud(T)$ is positive,
 then there exists a character 
$\chi:G_K\to F^*$ such that $\rho_2\simeq \rho_1\otimes
\chi$.

\item   Suppose that   $ ud(T) > 1/2$.   
Then there exists a finite extension $L$ of $K$
 and a character $\chi:G_L\to F^*$
such that $\rho_2|G_L\simeq \rho_1|G_L\otimes \chi$. 

\item Suppose $R=T^k: GL_2\to GL_{2^k}$. 
Assume that $ud(T)>1-2^{-(2m+1)}$.   that  the representations  $T^k(
\rho_1)$ and $T^k(\rho_2 )$ satisfy {\bf DH1}, for some positive
integer $k$ and  $\rho_1$ is
irreducible. Then there exists a character $\chi:G_K\to F^*$
such that $\rho_2\simeq \rho_1\otimes \chi$. 
\end{enumerate}
\end{corollary}
\begin{proof} Part 1) follows from Part (1) of Theorem
\ref{thm:special}. Part 2) follows from Theorem \ref{thm:smo} and
Part (1) of Theorem
\ref{thm:special}. We need to prove
Part (3) only for CM forms. This 
from Theorem \ref{thm:smo} and Proposition \ref{prop:rec-ind} (since
the representations are two dimensional, the subgroup $G'$ as in
Proposition \ref{prop:rec-ind} is of index $2$ in $G$, hence normal in
$G$).

\end{proof}

We now give the proof of Theorem \ref{thm:mf} stated in the
introduction. 
\begin{proof}[Proof of Theorem \ref{thm:mf}] For a newform $f$ of weight
  $k$, let 
\[\rho_f: G_{\Q}\to GL_2(\Q_l),\]
 be the
  $l$-adic representations of $G_{\Q}$ associated by the work of
  Shimura, Ihara and Deligne. 
   It has been shown by Ribet in \cite{Ri}, that the
  representation $\rho_f$ is semisimple, and the Zariski closure $G_f$
  of the image $\rho_f(G_{\Q})$ satisfies the following properties:
\begin{itemize}
\item If the weight of $f$ is at least two, then $G_f$ contains a
  maximal torus. 
\item If $f$ is a non-CM form, then $G_f=GL_2$.
\end{itemize}
Theorem \ref{thm:mf} follows by Corollary \ref{cor:twodim}
\end{proof}

\begin{remark}  A similar statement as below can be made for the class of
Hilbert modular forms too. A more general statement can be made when
the forms do not have complex multiplication, where we can take $R$ to
be an arbitrary rational representation of $GL_2$ with kernel
contained inside the centre of $GL_2$, rather than restrict ourselves
to symmetric and tensor powers. 
\end{remark}

We now consider an application to
Asai representations associated to holomorphic Hilbert modular
forms. With notation as in Example \ref{asai}, assume
further that $K/k$ is a quadratic extension of totally real number
fields, and that $\rho_1, ~\rho_2$ are $l$-adic representations
attached  respectively to holomorphic Hilbert modular forms $f_1, ~f_2$
over $K$. Let $\sigma$ denote a generator of the Galois group of
$K/k$, and further assume that $f_1$ (resp. $f_2$)  
is not isomorphic to any twist of
$f_1^{\sigma}$ (resp. $f_2^{\sigma}$ by a character, where
$f_1^{\sigma}$ denotes the form $f_1$ twisted by $\sigma$. Then the
twisted tensor  automorphic representations of $GL_2(\A_K)$, defined
by $f_i\otimes f_2^{\sigma}$ (with an abuse of notation) is
irreducible. Further it has been shown to be modular by
D. Ramakrishnan \cite{DR2}, and descend to define automorphic
representations denoted respectively by  
 $As(f_1)$ and $As(f_2)$ on $GL_2(\A_k)$. As a simple consequence of
Part d) of Theorem \ref{thm:special}, we have the following corollary:

\begin{corollary} Assume further that $f_1$ and $f_2$ are non-CM
forms, and that the Fourier coefficients of $As(f_1)$ and $As(f_2)$ are
equal at a positive density of places of $k$, which split in $K$. Then
there is an idele-class character $\chi$ of $K$ satisfying
$\chi\chi^{\sigma}=1$ such that,
\[ f_2=f_1\otimes \chi \quad \text{or}\quad  f_2 =f_1^{\sigma}\otimes
\chi.\]
\end{corollary}

\end{document}